\documentstyle[11pt]{article}
\input amssym.def
\input amssym
\input amssym.def
\input amssym

\def\R{\Bbb R}
\def\N{\Bbb N}

\newcommand{\fg}{\mbox{${\frak g}$}}

\newcommand{\fr}{\mbox{${\frak r}$}}

\newcommand{\ft}{\mbox{${\frak t}$}}

\def\la{\langle}
\def\ra{\rangle}

\def\tm{{\widetilde{M}}}
\def\tj{{\tilde{J}}}
\def\tom{{\widetilde{\om}}}
\def\tx{{\tilde{x}}}
\def\ty{{\tilde{y}}}
\newcommand{\q}{\mbox{${\rm\scriptscriptstyle q}$}}
\def\tX{{\widetilde{X}}}
\newcommand{\tq}{\mbox{$\tj^{\rm\scriptscriptstyle q}$}}
\newcommand{\tjj}{\mbox{$\tm_{\scriptscriptstyle \tj}$}}
\newcommand{\pjj}{\mbox{$\pi_{\scriptscriptstyle \tj}$}}

\newcommand{\C}{\mbox{${\Bbb C}$}}
\newcommand{\Z}{\mbox{${\Bbb Z}$}}
\newcommand{\p}{\mbox{${\rm p}$}}
\newcommand{\orb}{\mbox{${\rm\scriptscriptstyle orb}$}}
\newcommand{\Alb}{\mbox{${\rm Alb}$}}

\newcommand{\cf}{\mbox{${\rm cf}$}}
\newcommand{\ppi}{\mbox{${\pi^{\orb}_{1}(M)}$}}

\def\a{\alpha}
\def\G{\Gamma}

\def\lw{\longrightarrow}
\def \lo{\longmapsto}
 
 \def\om{{\omega}}

\def\s{\subset}

\newcommand{\pf}{\noindent{\it Proof.}\ }
\newcommand{\qed}{\hfill $\Box$ }
\def\sect#1{\section*{\centerline{\large\bf  #1}}}
\def\jielun#1{{\par\noindent\hskip 0.5 cm}
{\bf #1}\bgroup\it}
\def\endjielun{\egroup\par\bigskip}
\def\bdefi#1{{\par\noindent\hskip 0.5 cm}
{\bf #1}\bgroup\rm}
\def\edefi{\egroup\par\bigskip}
\def\beginli#1{{\par\noindent\hskip 0.5 cm}
{\bf #1}\bgroup\rm}
\def\endli{\egroup\par\bigskip}
\begin{document}
\title{\bf   {\CC {\char64}n K+ 4z
J{\char125} Sk A? WS 4z J{\char125} 5D JU Ku}}
\title{{\bf {Symplectic Convexity for Orbifolds}}
\thanks{Research partially supported by  SFB-237 of the DFG.}
\thanks{\it Mathematics Subject Classifications(2000): \ {\rm 53D05, 53D20}}}
\author{Qi-Lin Yang\\
Department of Mathematics,
Tsinghua University, \\
Beijing 100084, P. R. China\\
(E-mail: qlyang@math.tsinghua.edu.cn)\\
  }
\maketitle
\date{}

\begin{abstract}
We generalize symplectic convexity theorems for Hamiltonian
actions with proper momentum maps to symplectic actions on
orbifolds with mod-$\Gamma$ proper momentum maps.
\end{abstract}
\sect{\S1. \  Introduction}
An orbifold is a Hausdorff topological space
 locally modelled on $\R^{n}$ modulo finite
group actions. If the group actions are trivial, we recover the
concept of manifold. Quite an interesting thing is that the
enlarged category is closed under
 the quotients by finite groups.
 In symplectic geometry, an important construction of symplectic
quotients  called Marsden-Weinstein quotients,
 generically, are not
manifolds but symplectic orbifolds [1].
Naturally, we would like to generalize some
basic results on symplectic manifolds to the orbifold cases.

Atiyah, independently,  Guillemin and Sternberg established
symplectic convexity theorems  for Hamiltonian  torus actions on
symplectic manifolds
 in [2, 3, 4].
 Lerman and Tolman got
 the orbifold versions in [1].
In this note we give some generalizations of
their theorems using different methods.

\jielun{Theorem 1.1} Let $T$ be a torus and $(M,\om)$ a
 connected  symplectic $T$-orbifold. Let $\tm\lw M$ be
 the universal branch covering orbifold and
 $\G=\pi_{1}^{\orb}(M)$ the orbifold fundamental group.
Assume there is a momentum map $\tj:\tm\lw \ft^{*}$
 for the lifted action.
If $\tj$ is mod-$\G$ proper and the lifted $T$-action
 commutes with
that of $\G,$ then $\tj(\tm)$ is a closed convex set and
$\tj:\tm\lw \tj(\tm)$ is an open, fibre connected map.
\endjielun

Non-abelian version of Theorem 1.1
in K\"ahler and projective algebraic manifold cases
were independently proved by
Guillemin and Sternberg in [4]
 and Mumford in [5]. Kirwan [6] firstly
 accomplished the proof for non-abelian Hamilton
  action of a compact Lie group on an arbitrary connected compact symplectic manifold.
Sjammar [7], Heinzner-Huckleberry  [8] discussed
 extensions in algebraic and K\"ahler spaces.
 Flaschka-Ratiu [9] extended the results to the
  setting of {\it Poisson actions} of compact Poisson-Lie groups
 on symplectic manifolds.  If it is {\it symplectic action}
  on orbifold, we have the following extension:

\jielun{Theorem 1.2} Let $G$ be a connected compact Lie group and
$(M,\om)$ a connected  symplectic $G$-orbifold. Let $\tm\lw M$ be
the universal branch covering orbifold
 and $\G=\pi_{1}^{\orb}(M)$ the orbifold fundamental group.
Assume there is a $G$-equivariant momentum map
 $\tj:\tm\lw \fg^{*}.$
If $\tj$ is mod-$\G$ proper  and the lifted $G$-action commutes
with that of $\G,$  then $\tj(\tm)\cap \ft^{*}_{+}$ is a closed
convex  set and $\tj:\tm\lw \tj(\tm)$ is a fibre connected map.
\endjielun

There are several ways to show symplectic convexity theorems.
  Atiyah, Guillemin and sternberg,
 Lerman and Tolman , Kirwan in [2, 3, 4, 1, 6], employed Morse
theory.
 It is easy to show, using a normal form for Hamiltonian action,
the momentum map is locally convex [2, 3].  The Morse theory
gives rise to a global convexity theorem. Hilgert-Neeb-Plank [10]
offered  another proof by using a `local-global-principle',
dropped compactness of  the acted manifold  by an assumption
 that the momentum map is proper.
 Lerman-Meinrenken-Tolman-Woodward [11],
used the symplectic cutting technique.  It is a kind
of symplectic compactification. Intuitively, by
 cutting out infinity and collapsing the incision to a point we
 get a compact symplectic space such that the original
 non-compact symplectic space is equivariantly embedded in it as an
 open submanifold. Thus the proof is reduced to the compact case.
 This method work well in
orbifold cases and as a result the symplectic convexity theorems
are extended to non-compact orbifold cases [11].

In this paper, we use the techniques developed in [10] by
Hilgert-Neeb-Plank where the author dealt with  the manifold
cases. This proof is more analytical and elementary, it uses least
knowledge of symplectic geometry. In fact, we only need to know
that momentum map is locally convex, locally open, locally fiber
connected. But these are easily understood if we know the
symplectic version of slice theorem  for smooth groups actions.
Furthermore, this proof tells us clearly what causes the convexity
and why it should be so. To some extent, it builds the symplectic
convexity theorems on set theoretic topology.

Here is a brief description of the structure of this paper. In
Section 2 we review some basic concepts  and explain the
connections between symplectic action and Hamilton action.
Following the same idea of  Hilgert-Neeb-Plank, we define a
 map quotient $X_{f}$ for mod-$\G$ proper map
 in Section 3 and provethat it is a Hausdorff space. Finally we give  proofs of Theorem 1.1
and Theorem 1.2 in the last two sections.

\sect{\S2. Symplectic actions and Hamilton actions}
\hskip\parindent We refer to the Chapter 13 of [15] for a nice
account of  orbifolds and to [1] for  definitions of symplectic
orbifolds and Hamiltonian actions on them.  Let $G$ be a connected
Lie group with Lie algebra $\fg$ and $M$ a smooth connected
$G$-orbifold
 with symplectic structure $\om.$ A smooth action $G\times M\lw M$
is called {\it symplectic} if $\om$ is invariant under the action
of $G.$ In this case $M$ is called a {\it symplectic
$G$-orbifold}.   A symplectic action is called {\it Hamiltonian}
action if there exists a map $J:M\lw \fg^{*},$   called a {\it
momentum map},
 such that
$$i(\xi_{M})\om=dJ_{\xi},$$
where  $\xi_{M}$ is the infinitesimal generator
 corresponding to $\xi\in \fg$ and
 $J_{\xi}=\la J,\xi\ra$ denotes its $\xi$-component.
   A momentum map
  is called {\it equivariant} if it intertwines the
  action of $G$ on $M$
and the coadjoint action of $G$ on $\fg^{*}.$

A symplectic action is not always a Hamilton action. However, note
that $i(\xi_{M})\om$ is a closed form and
$i([\xi,\eta]_{M})\om=-d(\om(\xi_{M},\eta_{M})),$ so if
$H^{1}(M,\R)=0$ or $G$ is semi-simple, the symplectic action  is a
Hamilton action.

\jielun{Proposition 2.1} Let $(M,\om)$ be a connected symplectic
$G$-orbifold, and $G=R[G,G]$ a Levi-Malcev
 decomposition, here $R$ denotes the radical of $G.$ Then

{\rm (i)} there exists a $[G,G]$-equivariant momentum map $J:M\lw
\fg^{*};$

{\rm (ii)} the $G$-action is Hamiltonian if and only if the
$R$-action is Hamiltonian;

{\rm (iii)} if $M$ is a compact K\"ahler manifold with positive
Ricci curvature, in particularly, if $M$ is  Fano,
 then the $G$-action is Hamiltonian;

{\rm (iv)}if $M$ is a compact K\"ahler manifold,
 and  $G$ acts holomorphically symplectic on $M,$
then the $G$-action on $M$ is Hamiltonian if and only if the
 $R$-action on the Albanese variety
$\Alb(M) =H^{0}(M,\Omega^{1})^{*}/H_{1}(M,\Z)$ is trivial.
In particular, if $b_{1}(M)=0,$  the $G$-action is Hamiltonian.
\endjielun

\pf (i) and (ii) followed the discussions above; for (iii), note
in this case, $\pi_{1}(M)$ is a finite group, so we have
$H^{1}(M,\R)=0.$

For (iv),  first suppose that the $G$-action is Hamiltonian.
 Recall that the Albanese map $\a:M\lw\Alb(M)$ is equivariant.
To show the $R$-action on $\Alb(M)$ is trivial,
it suffices to show that every 1-parameter subgroup
 $r(t):=\la\exp(t\xi)\ra\s R$
 has a fix point on $M$  (\cf. [12, Proposition 1])£» here $\xi\in\fr$  and
 $\fr$ denotes the Lie algebra of $R.$ Let $J:M\lw \fr^{*}$ be the momentum map.
 Then the critical points of function $J_{\xi}=\la J,\xi\ra$,
which always exist since $M$ is compact, are the fixed points of
$r(t).$ Conversely, Suppose $R$ acts trivially on $\Alb(M),$ then
$R$ has fixed points in every fibre of $\a^{\scriptscriptstyle -1}
(\a(x))$ (cf. [13, Proposition]).  Thus  $\xi_{M}$
 has a zero point somewhere on $M$ for any $\xi\in \fr.$ Let
$\Xi=1/2(\xi_{M}-i{\Bbb J}\xi_{M})$ be the holomorphic vector
field defined by $\xi_{M},$  here ${\Bbb J}$ is the complex
structure of $M$. By (iii) of Theorem 1 in [14], there is a
function $f\in C^{\infty}(M,\C)$ such that
$i(\Xi)\om=\bar{\partial}f,$ so
 $i(\xi_{M})\om=i(\Xi)\om+\overline
 {i(\Xi)\om}=\bar{\partial} f+\partial{\bar{f}}.$
Let $f=\frac{1}{2}(g+ih), $  where $g,h\in C^{\infty}(M,\R).$ Then
$i(\xi_{M})\om=dg.$ So the $R$ action on $M$ is Hamiltonian. \qed

\vskip 0.5cm

In the following, let $M$ be an orbifold and $\p:\tm\lw M $ be the
universal branch cover. Then in general case $\tm$ is only an
orbifold  (\cf. [15, Chapter 13]). If $\tm$ is a manifold then $M$
is called a {\it good} orbifold. Let $\G:=\ppi$ be the orbifold
fundamental group of $M.$ Then $\G$ is a quotient group of
$\pi_{1}(M_{0})$ (\cf. [15, Chapter 13]),
 where $M_{0}$ is the regular points of $M.$
The action of $G$ lifted naturally on $\tm.$  Let
$\tom:=\p^{*}\om.$ Then $(\tm, \tom)$ is also
  a symplectic $G$-orbifold.  By Proposition 2.1
   there always exists a momentum map $\tj :\tm\lw\fg^{*}.$
  Still denote $\fr$ the radical
  of $\fg.$ Then $\fg=[\fg,\fg]\oplus \fr,$ the annihilator
  of $[\fg, \fg]$ in $\fg^{*}$ is $ [\fg,\fg]^{\circ}=\fr^{*}.$

\jielun{Proposition 2.2} Suppose that the actions of $G$ and $\G$
commute. Consider $\fg^{*}$  as a vector group,  then there exists
a homomorphism $h:\G\lw \fg^{*}$ such that
$$\tj(\gamma\cdot\tx)-\tj(\tx)=h(\gamma),\quad \forall \tx\in\tm.$$
If $\tj$ is $G$-equivariant, then $h(\tm)\subset
  \fr^{*}.$
 In particular, if $G$ is semisimple, then $\tj$ factors through
$\p$ so that there exists a $G$-equivariant momentum map $J:M\lw
\fg^{*}$
 such that $\tj =J\circ\p.$
\endjielun
\pf For any $\gamma\in \G$ and $\xi \in \fg$ set
$h_{\xi,\gamma}:=\tj_{\xi}\circ \gamma-\tj_{\xi}.$ If we denote
$\widetilde{H}_{f}$ the Hamiltonian vector field associated a
function $f$ on $\tm,$ then
$$\begin{array}{rcl}
\tom(\tx)(\widetilde{H}_{h_{\xi,\gamma}}(\tx),\tX(\tx)) &=&\la
d(\tj_{\xi}\circ\gamma)(\tx),\tX(\tx)\ra-
\la d\tj_{\xi}(\tx),\tX(\tx)\ra\\
&=&\la d\tj_{\xi}(\gamma \tx)\circ d\gamma(\tx),\tX(\tx)\ra-
\la d\tj_{\xi}(\tx),\tX(\tx)\ra\\
&=&\la d\tj_{\xi}(\gamma \tx),d\gamma(\tx)(\tX (\tx))
\ra-\la d\tj_{\xi}(\tx),\tX(\tx)\ra\\
&=&\tom(\tx)(\xi_{\tm}(\gamma\tx),d\gamma (\tX(\tx)))-
\tom(\tx)(\xi_{\tm}(\tx),\tX(\tx))\\
&=&(\gamma^{*}\tom)(\tx)(\xi_{\tm}(\tx),\tX(\tx))-
\tom(\tx)(\xi_{\tm}(\tx),\tX(\tx))\\
&=&0\\
\end{array}$$
since $\p\circ \gamma=\p$ implies
$\gamma^{*}\tom=\gamma^{*}(\p^{*}\om)= (\p\circ
\gamma)^{*}\om=\p^{*}\om=\tom.$  So $h_{\xi,\gamma}$ is
independent of $\tx$ and we can define $h(\gamma)\in\fg^{*}$ via
$\la h(\gamma),\xi\ra=h_{\xi,\gamma}.$ Clearly we have
$h(\gamma)=\tj\circ\gamma-\tj.$

For any $\tx\in\tm,$   clearly we have
$h(\gamma_{1}\gamma_{2}\tx)=\tj(\gamma_{1}\gamma_{2}\tx)-\tj
 (\tx)= (\tj(\gamma_{1}(\gamma_{2}\tx))-
\tj(\gamma_{2}\tx))+(\tj(\gamma_{2}\tx) -\tj(\tx)),
$
so  $h(\gamma_{1}\gamma_{2})=h(\gamma_{1})+h(\gamma_{2})$
and $h$ is a homomorphism.

If $\tj$ is $G$-equivariant, then for any $\tx\in \tm,$  we have
$Ad^{*}_{g}(h(\gamma))=Ad^{*}_{g}(h(\gamma\tx))=\tj(g\gamma\tx)-
\tj(g\tx)= \tj(\gamma g\tx)-\tj(g\tx)=h(\gamma g\tx)=h(\gamma),$
so $h(\tm)\subset {\fg^{*}}^{G}.$
 Since $\mu
\in {\fg^{*}}^{G}$ if and only if $ad^{*}_{\xi}\mu=0$   for any $
\xi\in \fg,$ that is $\la \mu,[\xi,\eta] \ra=0 $ for any $\xi$
 and $\eta\in\fg. $ Hence $h(\G)\s (\fg^{*})^{G}=[\fg,\fg]^{\circ}\cong
\fr^{*}.$

If $G$ is semisimple, then $\fr=0.$ Thus $\tj(\gamma
\tx)=\tj(\tx)$ for any  $\gamma\in\G.$ It follows that the
$G$-equivariant momentum map of $\tm$ descends to be a
$G$-equivariant map $J:M=\tm/\G\lw \fg^{*}.$ \qed

\sect{\S3. Quotient Space Modulo mod-$\G$ Map}
 \hskip\parindent
Let $X$ and $Y$ be topological spaces and  $f:X\lw Y$
 a continuous
map, $f$ is called {\it locally fibre connected} (\cf. [10,
Definition 3.1]) if for any $x\in Y$ there exist a neighborhood
$U$ of $x$ such that $f^{\scriptscriptstyle -1}(f(x))\cap U$ is
connected  for all $x\in U.$ If $f$ is locally fibre connected
map, the connected component of the fibre $f^{\scriptscriptstyle
-1}(f(x))$ passing through $x,$ denoted by $F_{x},$ is called a
{\it leaf} of $f.$ Define an equivalence relation $\sim$ on $X$ by
saying $x\sim y $ iff they belong to the same leaf. Let $X_{f}$
denote the quotient space with the quotient topology by shrinking
each leaf of $f$ to be a point. Then the quotient map, denoted by
$\pi_{f},$ is a continuous map. In general the structure of
$X_{f}$ is very complicate. For example, to assure $X_{f}$ be a
Hausdorff space, the {\it equivalence relation set} $E=:\{(x,y)\in
X\times X|x\sim y\}$ must be a closed subset. If $\pi_{f}$ is an
open map, $X_{f}$ is Hausdorff iff $E$ is closed. In [10], it is
proved that if $Y$ is a Euclid vector space and $f$ is a proper,
locally fibre connected and locally open, then
 $X_{f}$ is Hausdorff. In the following we will give a
 generalization of their result.

\bdefi{Definition 3.1} Let $H$ and $L$ be topological groups, $X$
a locally compact  topological $H$-space and $Y$ a locally compact
topological $L$-space. Let $f:X\lw Y$ be a continuous map and
$\rho: H\lw L$ be a continuous homomorphism. $f$ is called {\it
mod-$H$ proper} if for any  compact subset $C$ of $Y,$ there exist
a compact subset $B$ of $P$ such that $f^{\scriptscriptstyle
-1}(C)\s H\cdot B;$  and $f$ is called $\rho$-{\it equivariant,}
if $f(g\cdot x)=\rho(g)\cdot f(x)$ for any $g\in H$ and $x\in
X.$\edefi

\jielun{Proposition 3.1} Let $(X, d_{X})$ and $(Y, d_{Y})$ be
locally compact metric spaces. Assume that a locally compact group
$\G$ acts isometrically on $X$ and acts on $Y$ via an action
homomorphism $\rho:H\lw {\rm Iso} (Y).$ Suppose that $f:X\lw Y$ is
a locally fibre connected, locally open, mod-$\G$ proper,
$\rho$-equivariant, continuous map. Then
  $X_{f}$ is a Hausdorff topological space.
\endjielun

Clearly if $\G$ is a trivial group we recovered the result of
Hilgert-Neeb-Plank. To prove Proposition 3.1, we need the
following Lemma 3.1. For this we first give some notions that we
will use. For any closed subsets $A,B$ of $X,$ let
$d_{X}(x,B):=\inf_{y\in B}d_{X}(x,y)$ denote the distance from
$x\in X$ to $B$ and $d(A,B):=\sup_{x\in A}d_{X}(x,B).$ The {\it
Hausdorff distance} between $A$ and $B$ is defined by
$$d^{H}(A,B):=\max \{d(A,B),d(B,A)\}.$$
 Note that $d^{H}(A,B)=0$ if and only if $A=B.$
\jielun{Lemma 3.1}  Let $C$ be a compact  subset of $Y$ and
$W:=f^{\scriptscriptstyle -1}(C)=\G V,$ here  $V$ is a compact
subset of $X.$ Then

{\rm (i)} \ \ $\exists$  $ \delta > 0,$  $\forall x,y\in W,$  if
$d_{X} (x,y) \leq \delta$ and $f(x)=f(y),$ then $F_{x}=F_{y};$

{\rm (ii) }  $\exists$  $ \delta > 0,$  $ \forall \epsilon \in
[0,\delta],$ $ \exists$ $ \delta_{0} > 0,$ $\forall x,y\in W,$ if
$d_{X} (x,y)\leq \delta_{0}$ and $d_{Y}(f(x),f(y))\leq \delta,$
 then $d_{X}(x,F_{y})\leq \epsilon;$

{\rm (iii) }  $\forall$  $ \eta >0, $  $\exists $ $ \alpha >0,$ $
\forall x'\in V$ and $ y'\in W,$ if  $d_{X}(x',y')\leq\alpha,$
then $d_{Y}(f(x'),f(y'))\leq \eta;$

{\rm (iv) }  $\forall$  $\epsilon >0,$ $\exists$  $\beta>0,$
$\forall x,$ $ y \in W, $ if $d_{X}(x,y)\leq \beta,$ then
$d(F_{x},F_{y})\leq\epsilon.$
\endjielun
\pf Since  $f$ is $\rho$-equivariant, we have,
$$d_{Y}(f(\gamma \cdot x), f(\gamma \cdot y))=d_{Y}(
\rho(\gamma)\cdot f(x),\rho(\gamma)\cdot f(y))=d_{Y}(f(x),
f(y)).$$ Clearly $\gamma\cdot F_{x}=F_{\gamma \cdot x},$ so $
F_{x}=F_{y} $
 iff $\gamma\cdot F_{x}=\gamma\cdot F_{y}.$ These facts will be
 used throughout the following proof.

(i) We prove it by a contradiction. If otherwise,
 we can find series $\{x_{n}\},\{y_{n}\}$ in $W$
 such that $\lim_{n\rightarrow \infty}d_{X}(x_{n},y_{n})=0$
 and $f({x_{n}})=f({y_{n}}),$ but
$F_{x_{n}}\cap F_{y_{n}}=\emptyset.$
 Note that $f$ is mod-$\G$ proper,
we assume $x_{n}=\gamma_{n}x'_{n}$ such that $x'_{n}$ vary in a
compact subset and set $y_{n}=\gamma_{n}y'_{n}.$ Then
$f(x'_{n})=f(y'_{n})$ and $F_{x'_{n}}\cap F_{y'_{n}}=\emptyset.$

On the other hand, by extracting out subsequence if necessary, we
could assume $\lim _{n\rightarrow \infty}x'_{n}=x'_{\infty}\in W.$
Since $\lim_{n\rightarrow
\infty}d_{X}(x'_{n},y'_{n})=\lim_{n\rightarrow
\infty}d_{X}(x_{n},y_{n})=0,$ we get
 $\lim _{n\rightarrow \infty}y'_{n}=x'_{\infty}.$
So $\lim_{n\rightarrow \infty}F_{x'_{n}}\cap F_{y'_{n}}
=F_{x'_{\infty}}\not=\emptyset.$

(ii) If otherwise,  we can find $\epsilon_{0}>0$ and series
$\{x_{n}\},\{y_{n}\}$ such that for any  $\delta>0$ and
$d_{X}(x_{n},y_{n})\leq \delta$ and   $
d_{Y}(f(x_{n}),f(y_{n}))\rightarrow 0,$
 but $d_{X}(x_{n},F_{y_{n}})\geq \epsilon_{0}.$
Similar to the proof of (i), without loss of generality,  we may
assume $x_{n}\rightarrow x_{\infty}$ and $ y_{n}\rightarrow
y_{\infty},$
 so $d_{X}(x_{\infty},y_{\infty})\leq \delta$ and
$f(x_{\infty})=f(y_{\infty}).$ By (i),  if $\delta$ is small
enough then we have $F_{x_{\infty}}=F_{y_{\infty}}.$   Which means
$x_{\infty}$ and $y_{\infty}$ lie in the same leaf of $f.$
 So $\lim_{n\rightarrow{\infty}}d_{M}(x_{n},F_{y_{\infty}})=0.$
 A contradiction.

(iii) Since $f$ is continuous, by definition,
 $\forall\eta>0 $ and $\forall x'\in V,$
   there exists $\alpha_{x'}>0$ such that
 $\forall y'\in W$
if $d_{X}(x',y')\leq\alpha_{x'}$ then
$d_{Y}(f(x'),f(y'))\leq\eta.$ Let $\Sigma=\bigcup_{x'\in Y}
\{y'\in W|d_{X}(x',y')\leq\alpha_{x'}\}\cap V$ be an open
 cover of $V.$ Note $V$ is compact,
by Heine-Borel Covering Theorem there is a finite subcover
$\Sigma'=\bigcup_{x'_{i}\in V} \{y'\in
W|d_{X}(x'_{i},y')\leq\alpha_{i},i=1,\cdots ,n\}\cap V=V.$ Set
$\alpha=\min\{\alpha_{i}|i=1,\cdots n\}.$ So for any $x'\in Y$
 and $y'\in W$ if $d_{Y}(x',y')\leq\alpha,$ we have
$d_{Y}(f(x'),f(y'))\leq\eta.$

(iv) Let $\delta$ and  $\delta_{0}$ as in (ii). By (iii) there
exists $\beta>0,$  if $d_{X}(x,y)\leq \beta$ then
 $d_{Y}(f(x), f(y))\leq \delta.$
Now let $\alpha=\min\{\delta_{0},\beta\}$ and  $\epsilon\in
[0,\alpha/2],$ set $E_{x}:=\{z\in F_{x}|d_{X}(z,F_{y})\leq
\epsilon\}.$ Then $E_{x}$ is non-empty closed subset of $F_{x}.$
If we can show $E_{x}$ is open in $F_{x},$ then (iv) follows.

In fact, let $z\in E_{x}$ and $w\in F_{x}\cap\{w\in
W|d_{X}(w,z)\leq \epsilon\}.$ Let $w'\in F_{y}$ such that
$d_{X}(z,w')\leq \epsilon.$ Then $d_{X}(w,w')\leq 2 \epsilon \leq
\alpha$ and $ d_{Y}(f(w),f(w'))\leq\delta. $ By (ii) we have
$d_{Y}(w,F_{y})\leq \epsilon.$ Thus $E_{x}$ is open in $F_{x}.$
\qed

\vskip 0.8cm
 \noindent{\it Proof of Proposition 3.1.}  Since $f$ is
locally open, $\pi_{f}$ is an open map.
 It suffices to show  the equivalence relation set is closed, that
 is to show  for any series $\{x_{n}\},\{y_{n}\}$ in $X,$
 if $x_{n}\rightarrow x_{\infty} $
and $y_{n}\rightarrow y_{\infty}$ and $F_{x_{n}}=F_{y_{n}},$ we
have $F_{x_{\infty}}=F_{y_{\infty}}.$ In fact let
$C=\{f(x_{n}),f(y_{n})|n=1,2,\cdots,\infty\},$ then $C$ is a
compact subset of $Y.$ Let $W:=f^{\scriptscriptstyle -1}(C).$
Using (iv) of Lemma 3.1, we have $0\leq
d^{H}(F_{x_{\infty}},F_{y_{\infty}})\leq
 d^{H}(F_{x_{\infty}},F_{x_{n}})+
d^{H}(F_{x_{n}},F_{y_{n}})
 +d^{H}(F_{y_{n}},F_{y_{\infty}})
\rightarrow 0$ as $n\rightarrow \infty.$
 So $F_{x_{\infty}}=F_{y_{\infty}}.$  Thus $X_{f}$ is a Hausdorff topological space. \qed

\sect{\S4. Abelian Convexity} From now on we will continue the
discussions in Section 2. In the following, we use the same
notations and terminologies as in Section 3 unless otherwise is
especially stressed.  Let $M$ be a symplectic orbifold and $\tm$
its universal branch covering orbifold
 and $\G=\pi_{1}^{\orb}(M)$ the orbifold fundamental group.
  Assume $G=T$ is a torus, and $T\times M\lw M$ is a symplectic
   action and the lift $T\times \tm\lw \tm$ is a Hamilton
    action with momentum map $\tj.$  Then $\tj$ is locally
     fibre connected, locally open, continuous map. For these
    properties we refer to [1, 3, 10] for detailed accounts.
     We suppose that $\tj$ is mod-$\G$ proper. Choose a
      Riemannian metric on $M,$ we may lift it to $\tm$
      and assume $\G$ acts isometrically on $\tm$ relative
      to the lifted metric.
Let $\pjj:M\lw \tjj$ be the quotient map. $\tj:\tm\lw \ft^{*}$
induces a map $\tq:\tjj\lw \ft^{*}$ such that $\tq\circ \pjj=\tj.$

\jielun{Proposition 4.1}   $\tjj$ is a Hausdorff topological
space.
 \endjielun
\pf  We take $\ft^{*}$ as a metric space with Euclid metric. By
Proposition 2.2, $\tj\circ\gamma=\tj+h(\gamma).$ If we consider
$\G$ as an isometry transformation group acting on $\ft^{*}$ by
translation via homomorphism $h,$ then $\tj$ is a $h$-equivariant
map. Thus we can use Proposition 3.1 to conclude $\tjj$ is a
Hausdorff topological space. \qed

Any $\gamma\in \G$ descends to be a homeomorphism of $\tjj,$
denoted by $\tilde{\gamma},$ satisfying
$\tilde{\gamma}\circ\pjj=\pjj\circ\gamma.$ Let
$\widetilde{\G}:=\{\tilde{\gamma}|\gamma\in\G\}.$ Then $\tjj$
 is Hausdorff topological space by Proposition 3.1,
 so $\tq:\tjj\lw \ft^{*}$ is mod-$\widetilde{\G}$
 proper continuous map. Moreover  $\tq
  \circ{\tilde{\gamma}}-\tq$
  is a constant function on $\tjj$
for any  $\tilde{\gamma}\in \widetilde{\G}$  by Proposition 2.2.

A continuous map $c:[0,1]\lw \tjj $ is called a {\it regular
curve}  (\cf. [10, Definition 3.6]) if $\tq\circ c$ is piecewise
differentiable. The length of $\tq\circ c$ is denoted by
$l(\tq\circ c).$ For any $\tx^{\q}_{0},\tx^{\q}_{1}\in \tjj,$ let
\begin{center}
$d(\tx^{\q}_{0},\tx^{\q}_{1}):=\inf\{l(\tq\circ c)|c$  is a
regular curve,  $c(i)= \tx^{\q}_{i},i=0,1\}.$
\end{center}
Clearly $d$ is symmetric and satisfies the triangle inequality,
and $d_{\ft^{*}}(\tq(\tx^{\q}_{0}), \tq(\tx^{\q}_{1}))\leq
d(\tx^{\q}_{0},\tx^{\q}_{1}).$

 \jielun{Proposition 4.2}   For any
$\tx^{\q}\in \tjj$ and $r>0,$ the closed ball
$B_{r}(\tx^{q}):=\{{\tilde{y}}^{\q}\in\tjj|d({\tilde{y}}^{\q},
\tx^{\q})\leq r\}$ is compact.
\endjielun
\pf For any
$\tx^{\q}\in \tjj$ and $r_{\scriptscriptstyle 0}>0,$ let $B=B_{r_{0}}(\tq(\tx^{\q}))$
 be a closed ball in
 $\ft^{*},$ then there exist a compact subset $A$ of $\tjj$
  such that $(\tj^{\q})
^{\scriptscriptstyle -1}(B)\s
 {\widetilde{\G}}\cdot A.$ Let $A_{0}$ denote the connected
component of ${\widetilde{\G}}\cdot  A$ containing $\tx^{\q},$
then $A_{0}$ is a compact neighborhood of $\tx^{\q}.$
 So $\tjj$ is locally compact.
 Furthermore, $\exists $ $\delta>0$ such that
$B_{\delta}({\tilde{x}}^{\q})\s A_{0},$
so  $B_{\delta}({\tilde{x}}^{\q})$ is compact.
 We can find  finite such closed ball
 $B_{\delta}({\tilde{x}}^{\q})$ to cover
  $B_{r}(\tx^{\q}),$ hence   $B_{r}(\tx^{\q})$ is compact.\qed

\jielun{Proposition 4.3}   $d:\tjj\times \tjj\lw \R$ is a metric.
\endjielun
\pf It suffices to show $d$ separate points. We assume that $
d(\tx^{\q},\ty^{\q})=0.$ Then $\mu:=\tq(\tx^{\q})=\tq(\ty^{\q}).$
Let $B:=\tq^{\scriptscriptstyle -1}(\mu)$ and
$C:=\tj^{\scriptscriptstyle -1}(\mu)=\pjj^{\scriptscriptstyle
-1}(B)=\{F_{x_{i}}|F_{x_{i}}\cap F_{x_{j}}=\phi$ if $i\not=j$
$\}.$ We claim that $\{x_{i}\}$ has no convergence point.
Otherwise $x_{n}\rightarrow x_{\infty},$ by (i) of Proposition
3.1,  $\exists N>0$ such that when $n>N$ all $F_{x_{n}}$ coincide,
a contradiction. Thus we can find disjoint open sets
$\{U^{\q}_{i}\}$  in $\tjj,$  such that each $U^{\q}_{i}$ contains
only one element $\tx^{\q}_{i}$ of $\tq^{\scriptscriptstyle
-1}(\mu).$ Since $\tj$ is a locally convex map, clearly so is
$\tq.$ So $\tq(U^{\q}_{i})$ contain a closed ball with positive
radius $\epsilon$ and any regular curve $c$ starts at $x^{\q}_{i}$
and leaves $U^{\q}_{i}$ satisfies $l(\tq\circ c)\geq \epsilon.$ So
$\tx^{\q},\ty^{\q}$ must lie in the same closed ball. This shows
$\tx^{\q}=\ty^{\q}.$\qed

\jielun{Remark 4.1}
 {\rm  Proposition 4.2 and 4.3 together say $\tm_{\tj}$ is a connected
 locally compact metric space. Note that a metric space is not necessary
 locally compact. The simplest example is the rational
 number ${\Bbb Q}$ as a subspace of ${\Bbb R}^{1}$ with
 Euclid metric, clearly it is not locally compact since
 any compact subset of ${\Bbb Q}$ is a finite set. A connected
 Hausdorff space is not necessary locally compact, too. For example,
 the quotient of  ${\Bbb R}^{1}$ modulo the equivalence
  relation set $E={\Bbb Z}\times {\Bbb Z}$ is clearly connected
  but
   not locally compact.}
\endjielun
\vskip 0.8cm

\noindent{\it Proof of Theorem 1.1.} Fix $\tx^{\q},\ty^{\q}\in
\tjj$ and let $d_{0}:=d(\tx^{\q},\ty^{\q}).$ For any $n\in \N,$
there exist a regular curve $c_{n}$ connecting $\tx^{\q}_{0}$ and
$\tx^{\q}_{1}$ with $l(\tq\circ c_{n})\leq d_{0}+\frac{1}{n}.$ Let
$\tx^{\rm\scriptscriptstyle q;n}_{\rm\scriptscriptstyle 1/2}$ be
the midpoints of $c_{n}.$ They are contained in the ball
$B_{2d_{0}}(\tx^{\q})$ which is compact, so they have a coherent
point $\tx^{\q}_{\rm\scriptscriptstyle 1/2}.$ This point satisfies
$$d(\tx^{\q},\tx^{\q}_{\rm\scriptscriptstyle 1/2})
=d(\tx^{\q}_{\rm\scriptscriptstyle 1/2},\ty^{\q})=d_{0}/2.$$
Repeat this process for the pairs of points
$(\tx^{\q},\tx^{\q}_{\rm\scriptscriptstyle 1/2})$ and
$(\tx^{\q}_{\rm\scriptscriptstyle 1/2},\ty^{\q})$ to obtain
$\tx^{\q}_{\rm\scriptscriptstyle 1/4}$
 and $\tx^{\q}_{\rm\scriptscriptstyle 3/2}$ respectively, satisfying

$$d(\tx^{\q},\tx^{\q}_{\rm\scriptscriptstyle 1/4})
=d(\tx^{\q}_{\rm\scriptscriptstyle 1/4},\tx^{\q}_{\rm\scriptscriptstyle 1/2})
=d(\tx^{\q}_{\rm\scriptscriptstyle 1/2},
\tx^{\q}_{\rm\scriptscriptstyle
3/4})=d(\tx^{\q}_{\rm\scriptscriptstyle 3/4},\ty^{\q}).$$
Inductively we find points $\tx^{\q}_{\rm\scriptscriptstyle
k/2^{m}},$  for $0\leq k\leq 2^{m}$ such that
$$d(\tx^{\q}_{\rm\scriptscriptstyle k/2^{m}},
\tx^{\q}_{\rm\scriptscriptstyle k'/2^{m}})=d_{0}|k/2^{m}-k'/2^{m}|.$$
So we can extend $k/2^{m}\lo \tx^{\q}_{\rm\scriptscriptstyle
k/2^{m}}$ to a continuous map $c:[0,1]\lw \tjj$ such that
$d(c(t),c(t'))=d_{0}|t-t'|.$ This means
$$d_{\ft^{*}}(\tq\circ c(t),\tq\circ c(t'))=d_{0}|t-t'|$$
which can  only happen iff $\tq\circ c$ is a straight line. So
$\tj(\tm)=\tq(\tjj)$ is a convex set.

To show that the fibres of $\tj$ are connected, it suffices to
show $\tq$ is injective.  We assume $\tq(\tx^{\q})=\tq(\ty^{\q}).$
We construct a regular curve $c$ connecting $\tx^{\q}$ and
$\ty^{\q}$ as in the previous paragraph. Then
$d(\tx^{\q},\ty^{\q})=d_{\ft^{*}}(\tq\circ c(0), \tq\circ
c(1))=0.$ So that $\tx^{\q}=\ty^{\q}.$ In view of what we have
already shown, $\tq$ is a homeomorphism. Since $\pi_{\tj}$ is an
open map, $\tj=\tq\circ \pi_{\tj}$ is an open map as well.  Thus
$\tj:\tm\lw \tj(\tm)$ is an open, fibre connected map.

For any $\mu\in \overline{\tq(\tm)},$ assume
$\lim_{n\rightarrow +\infty}\tq(\tx^{\q}_{n})=\mu.$ Since
$\tq$ is a mod-$\widetilde{\G}$
 proper  map, we can assume
 $\tx^{\q}_{n}=\tilde{\gamma}\cdot\ty^{\q}_{n}$ such that
$\ty^{\q}_{n}$ varies in a compact subset.
By extracting out subsequence if necessary, we
could assume $\lim _{n\rightarrow \infty}\ty^{\q}_{n}=
\ty^{\q}_{\infty}\in\tm.$ Thus
$\mu=\tq(\ty^{\q}_{\infty})+\lim _{n\rightarrow
\infty}h(\tilde{\gamma_{n}})\in \tq(\tm)+h(\widetilde{\G})=\tq(\tm).$
Hence $\tj(\tm)$ is a closed convex set.
 \qed

\sect{\S5. Non-Abelian Convexity} First we review the symplectic
cross-section theorem for actions of compact connected Lie group
$G.$  (\cf. [3, Theorem 6.4] and [11, Theorem 3.1]).  Let
$T\subset G$ be a maximal torus and $\ft^{*}_{+}\subset \ft^{*}$ a
positive Weyl chamber. For any $\lambda\in \fg^{*},$
 there is a unique point in $\ft^{*}_{+}$ which is the
intersection point of the coadjoint orbit $Ad^{*}_{G}\lambda$ and
$\ft^{*}_{+}.$ Thus $\ft^{*}_{+}$ parameterizes the coadjoint
orbits and  is a {\it  section} for the coadjoint action. Now let
$M$ be a connected Hamilton $G$-orbifold with equivariant moment
map $J:M\lw \fg^{*},$ we can ``pull back" the  section for
coadjoint action to a section for the $G$-action on $M$ via  $J.$
Let $\sigma$ denote the interior of $\ft^{*}_{+}.$
 The preimage $Y:=J^{\scriptscriptstyle -1}(\sigma)$
is a connected $T$-invariant {\it suborbifold} of $M.$ The {\it
Symplectic Section Theorem} claim that $Y$ is a {\it symplectic}
suborbifold, thus a ``symplectic section"(\cf. [11, Theorem 3.1]).
It is easy to see the restriction $J_{Y}$ of $J$  to $Y$ is a
moment map for action of  $T$ and $G\cdot Y$ is dense in $M.$ Thus
the ``symplectic section" set up a bridge between
  torus action and non-abelian group action. We will
  use these facts to prove  Theorem 1.2 via using Theorem 1.1.

\vskip 0.8cm \noindent{\it Proof of Theorem 1.2} Let $\sigma$ be
the interior of the Weyl chamber $\ft^{*}_{+}$ and $\widetilde{Y}=
:\tj^{\scriptscriptstyle -1}(\sigma)$ the symplectic section.
$\widetilde{Y}$ is a symplectic $T$-orbifold with momentum map
$\tj_{\widetilde{Y}}.$
 Since $\sigma$ is a relative  open  subset of $\ft^{*}_{+},$
   we can choose an ascending sequence of closed subsets
   $\sigma_{i}\s \sigma$ such that
$\cup_{i\in \N}\sigma_{i}=\sigma.$ Let $\widetilde{Y}_{i}:=
\tj^{\scriptscriptstyle -1}_{\widetilde{Y}}(\sigma_{i})$ be the
closed subsets of $\tjj.$ Since $\widetilde {Y}$ is connected, we
can choose an ascending sequence of connected components
$\widetilde {Y}'_{i}$ of $\widetilde {Y}_{i}$ such that
$\widetilde {Y}=\cup_{i\in\N} \widetilde{Y}'_{i}.$ The restriction
$J|_{\widetilde {Y}'_{i}}:\widetilde {Y}'_{i}\lw \sigma_{i}$ is a
mod-$\widetilde{\G}$-proper. Clearly $J|_{\widetilde {Y}'_{i}}$ is
also a locally fibre connected,
 locally convex, locally open map.
By Theorem 1.1, we know $\tj(\widetilde {Y}'_{i})$  form an
   ascending sequence of closed convex subsets of $\sigma.$
   Hence $\tj(\widetilde {Y})$ is  convex.
  So $\overline{\tj(\tm)}\cap
   \ft^{*}_{+}=\overline{\tj(Y)}$ is convex locally polyhedral
    set. If we can prove $\tj(\tm)$ is closed, then
    ${\tj(\tm)}\cap \ft^{*}_{+}=\overline{\tj(Y)}$ is a
    closed convex set. In fact, since $\tj$ is mod-$\G$
proper, the quotient map $\widehat{J}:M=\tm/\G\lw \fg^{*}/h(\G)$
is proper, $\widehat{J}(M)$ is closed in $\fg^{*}/h(\G),$
 so $\tj (\tm)$ is a closed subset of $\fg^{*}.$

 It remains to show that the fibre $\tj^{\scriptscriptstyle -1}
 (\mu)$  is connected for any $\mu\in\fg^{*}.$ By Theorem 1.1,
 the fibres of $\tj_{\widetilde{Y}}$  are connected.
  Since $G\cdot Y$ is dense in $M$ and $\tj$ is equivariant,
  $\tj|_{\scriptscriptstyle G\cdot Y}$  is fibre connected.
Clearly $\tj^{\scriptscriptstyle -1}( Ad^{*}_{G}\mu)=G\cdot \tj^
{\scriptscriptstyle -1}(\mu).$ Since  $G$ and  $G_{\mu}$ are
connected, $\tj^ {\scriptscriptstyle -1}(\mu)$ is connected if and
only if $\tj^{\scriptscriptstyle -1}( Ad^{*}_{G}\mu)$ is
connected. So we may assume $\mu\in \ft^{*}_{+}.$

Now for any $\mu\in \sigma,$ the fibre
 $\tj^{\scriptscriptstyle -1}(\mu)$ is connected by Theorem 1.1.  So for any
 convex open neighborhood $B$ of $\mu,$ the set
 $G\cdot (\tj^{\scriptscriptstyle -1}(B)\cap \sigma)$
 is connected. Since $\tj^{\scriptscriptstyle -1}(Ad^{*}_{G}
\cdot B\cap \ft^{*}_{+})\cap G\cdot \tj^{\scriptscriptstyle -1}
(\sigma)=G\cdot \tj^{\scriptscriptstyle -1}(B\cap \sigma),$ we
know $\overline{\tj^{\scriptscriptstyle -1} ( Ad^{*}_{G}\cdot
B\cap \ft^{*}_{+})}=\overline{ G\cdot \tj^{\scriptscriptstyle
-1}(B\cap \sigma)}$ is connected.

For any $\mu\in \ft^{*}_{+},$
let $B_{i}$ be  convex open Neighborhoods  with
$\mu\in \overline{B_{i}}$ and $B_{i+1} \s B_{i}$
 such that $\cap_{i\in \N}
\overline{B_{i}}=\{\mu\}.$
Then $\tj^{\scriptscriptstyle -1}( Ad^{*}_{G}\mu)=
\cap_{i\in \N}\overline{\tj^{\scriptscriptstyle -1}
( Ad^{*}_{G}\cdot B_{i}\cap \ft^{*}_{+})}$ is connected.
It  follows that $\tj^{\scriptscriptstyle -1}(\mu)$ is
connected. So $\tj$ is a  fibre connected map.
\qed
\section*{Acknowledgement}

 Part of this work was done while I was
 visiting as a guest fellow at
 the Institut f\"ur Mathematik, Ruhr Universit\"at Bochum, Germany.
I would like to thank Prof. A. Huckleberry for encouragements
and P. Heinzner  for  useful discussions.
I also wish to thank the referee for many skillful
 comments and for pointing
out an error in a previous version of the manuscript.

\begin{center}
{\bf\Large References}
\end{center}
{\parindent=0pt

\def\toto#1#2{\centerline{\makebox[1.5cm][l] {#1\hss}
\parbox[t]{13cm}{#2}}\vspace{\baselineskip}}

\toto{[1]}{{\rm Lerman, E. \& Tolman, S.,} {Hamilton Torus actions
on Symplectic Orbifolds and Toric Varieties,} {\it Tran. A. M.
S.,} {\bf 349} (1997), 4201-4230.}

\toto{[2]}{{\rm Atiyah, M.,}
{Convexity and commuting Hamiltonians,}
 {\it Bull. Lond. Math. Soc.,} {\bf 14} (1982), 1-21.}

\toto{[3]}{{\rm  Guillemin, V. \&   Sternberg, S.,}
 { Convexity properties of moment mapping,} {\bf I,}
{\it Inv. Math.,} {\bf 67} (1982), 491-513.}

   \toto{[4]}{{\rm  Guillemin, V. \&   Sternberg, S.,}
 { Convexity properties of moment mapping,} {\bf  II,}
{\it Inv. Math.,} {\bf 77} (1984), 533-546.}

  \toto{[5]}{{\rm  Ness, L.,}
{A stratification of the null cone via the moment map,
with an appendix by Mumford, D.,}
 {\it Amer. J. Math.,} {\bf 106} (1984), 1281-1330.}

\toto{[6]}{{\rm Kirwan, F.,}
{ Convexity properties of moment mapping,} {\bf  III,}
{\it Inv. Math.,} {\bf 77} (1984), 547-552.}

\toto{[7]}{{\rm  Sjamaar, R.,}
  { Convexity properties of moment mapping re-examined,}
{\it Adv. Math.,} {\bf 138} (1998), 46-91.}

\toto{[8]}{ {\rm Heinzner, P. \& Huckleberry, A.,} {K\"ahler
potentials and convexity properties of the moment map,} {\it Inv.
Math.,} {\bf 126} (1996), 65-84.}

\toto{[9]}{ {\rm Flaschka, H. \& Ratiu, T.,} {A convexity theorem
for Poisson actions of compact Lie groups,} {\it Ann. Sci.
l'\'Ecol Norm. Sup\'er.,} {\bf 29} (1996), 787-809.}

\toto{[10]}{ {\rm Hilgert, J. \& Neeb, K,-H. \& Plank, W.,}
{Symplectic Convexity Theorems and Coadjoint Orbits,}
{\it Compo. Math.,} {\bf 94} (1994), 129-180.}

\toto{[11]}{{\rm Lerman, E. \& Meinrenken, E. \& Tolman, S. \&
Woodward, C.,} {Non-abelian convexity by symplectic cuts,} {\it
Topo.} {\bf 37} (1988), 245-259.}

\toto{[12]}{ {\rm  Huckleberry, A.\& Wurzbacher, T.,}
{Multiplicity-free complex manifolds,} {\it Math. Ann.,}{\bf 286}
(1990), 261-280.}

\toto{[13]}{{\rm  Sommese, A. J.,} {Extension theorems for
reductive group actions on compact K\"ahler Manifolds,} {\it Math.
Ann.,} {\bf 218} (1975), 107-116.}

\toto{[14]}{{\rm LeBrun, C. \& Simanca, S. R.,} {Extremal K\"ahler
Metrics and Complex Deformation Theory,} {\it Geom. Funct. Anal.,}
{\bf 4} (1993), 298-336.}

\toto{[15]}{{\rm Thurston, W.,}
 {The Geometry and Topology of 3-manifolds,}
 {Mimeographed notes} {Princeton University.}}
}

\newpage

{\Huge \begin{center}
     {\rm To}:\\
     {\bf Prof.\ Dr.\ Victor\ Bangert}\\
     {\sf Mathematisches \  \ Institut}\\
     {\sf Universität\ \  Freiburg}\\
     {\sf Eckerstrasse \ \ 1,\ \ Raum 336}\\
     {\sf 79104\ \ Freiburg \ \ im\ \ Breisgau}\\
     {\sf Germany}\\
    \end{center}

}

\end{document}